\documentclass{article} 

\usepackage[english]{babel} 
\usepackage{amssymb}
\usepackage{amsmath}
\usepackage{txfonts}
\usepackage{mathdots}
\usepackage[classicReIm]{kpfonts}
\usepackage{graphicx}
\usepackage{color}   

\def\e{\hbox{e}}
\def\exp{\hbox{exp}}
\def\ds{\displaystyle}

\def\CC{\mathbb{C}}

\def\erf{\hbox{erf}}

\def\L{\mathcal L} 

\def\e{\hbox{e}}

\def\LL{\mathcal L} 
\def\ds{\displaystyle}

\newcommand{\bee}{\begin{equation}}
\newcommand{\ee}{\end{equation}}

\begin{document}




\title{Application of the Efros theorem to the function represented by the inverse Laplace transform of $ s^{-\mu}\, \exp(-s^\nu)$}
\date{22 February 2021}
\maketitle
%
\noindent \textbf{\Large
Alexander Apelblat${}^{a}$ and 
Francesco Mainardi${}^{b \ast}$}

\noindent ${}^{a}$Department of Chemical Engineering, Ben Gurion University of the Negev, 84105 Beer Sheva, 84105, Israel. Email: apelblat@bgu.ac.il

\noindent ${}^{b}$Department of Physics and Astronomy, University of Bologna, Via Irnerio 46, 40126 Bologna, Italy. Email: francesco.mainardi@bo.infn.it  \\
ORCID:  0000-0003-4858-7309

\noindent
$^\ast$ Corresponding author. 

\noindent 
\\
\noindent \textbf{ Keywords: }Efros theorem; inverse Laplace  transforms; Wright functions; Mittag-Leffler functions; Volterra functions; modified Bessel functions;  finite, infinite and convolution integrals.\\ \\
{\bf This paper has been published in the  Special
Issue of Symmetry (MDPI)  
\\ "Special Functions and Polynomials", with guest Editor  Paolo Emilio Ricci,
\\ Vol 13 (2021), art 354, 15 pages.
DOI: 10.3390/sym13020354}

\noindent 
\\ 
\noindent \textbf{ AMS Subject Classification : }26A33, 33C10, 33E12, 34A25, 44A20 

\noindent \textbf{}
\\
\noindent \textbf{\large Abstract:} 
Using a special case of the Efros theorem which was derived by Wlodarski, and operational calculus, it was possible to derive many infinite integrals, finite integrals and integral identities for the function represented by the inverse Laplace transform. 
The integral identities are mainly in terms of convolution integrals with the Mittag-Leffler and Volterra functions. 
The integrands of determined integrals include elementary functions (power, exponential, logarithmic, trigonometric and hyperbolic functions) and the error functions, the Mittag-Leffler functions and the Volterra functions. 
Some properties  of the inverse Laplace transform of  
$s^{-\mu} \exp(-s^\nu)$ with $\mu \ge0$ and $0<\nu<1$
 are presented.

%
%
\section{Introduction}
\noindent Inversions of the Laplace transforms of exponential functions
\begin{equation} \label{GrindEQ__1_} 
\begin{array}{l} 
{f_{\nu ,\mu } (t)=  \LL^{-\, 1} \left\{F(s)\right\}
=\dfrac{1}{2 \pi  i} \, {\ds \int _{c\, -\, i \infty }^{c\, +\, i \infty }}\e^{s t}\, F(s)\, ds
\quad  ;\quad c>0,}
 \\   \phantom{\rule{1pt}{15pt}}  
 {F(s)=\dfrac{\e^{- s^{\nu } } }{s^{\mu } } \quad ;\quad 0<\nu <1\quad ;\quad \mu \ge 0,} \end{array} 
\end{equation} 
were during the 1945 - 1970 period in a focus of attention of a number of well-known mathematicians like Humbert [1], Pollard [2], Wlodarski [3], Mikusinski [4-7], Wintner [8], Ragab [9] and Stankovi\v{c} [10]. In the case $\mu=0$, Mikusinski was able to obtain the inverse Laplace transform in terms of integral representations [7]
\begin{equation} \label{GrindEQ__2_} 
\begin{array}{l} 
{\LL\left\{f_{\nu ,0} (t)\right\}={\ds \int _{0}^{\infty} }\e^{-\, s t}\, f_{\nu ,0} (t)\, dt
=\e^{-\, s^{\nu } }  \quad ;\quad 0<\nu <1\quad ;\quad t>0} 
\\   \phantom{\rule{1pt}{15pt}}   
{f_{\nu ,0} (t)=
\dfrac{1}{\pi } \,{\ds  \int _{0}^{\infty} }\e^{-\, u t}  \e^{-\, u^{\nu } \cos (\pi  \nu )} \sin [u^{\nu } \sin (\pi  \nu )]\,  du} 
\\   \phantom{\rule{1pt}{15pt}}  
{f_{\nu ,0} (t)
=\dfrac{2}{\pi } \, {\ds \int _{0}^{\infty }}\e^{-\, u^{\nu } \cos (\frac{\pi  \nu }{2} )}  \cos [u^{\nu } \cos (\frac{\pi  \nu }{2} )]\,  \cos (u t)\, du} 
\end{array} 
\end{equation} 
and also as the finite trigonometric integral
\begin{equation} \label{GrindEQ__3_} 
\begin{array}{l} 
{f_{\nu ,0} (t)
=\dfrac{\nu }{\pi  (1-\nu )  t} \, {\ds \int _{0}^{\pi} }\xi  \e^{-\, \xi }\,  du\quad ;\quad 
0<\nu <1\quad ;\quad t>0} 
\\ \phantom{\rule{1pt}{15pt}}  
{\xi =\dfrac{1}{t^{\nu /(1\, -\, \nu )} } \, 
\left(\dfrac{\sin (\nu  u)}{\sin u} \right)^{\nu /(1\, -\, \nu )} \dfrac{\sin [(1-\nu ) u]}{\sin u} } \end{array} 
\end{equation} 
It was established that the functions $f_{\nu,0}(t)$
can be expressed in terms of exponential and parabolic cylinder functions when 
$\nu  = 1/2$  [9,11]
 and by help of  the Airy functions and their first derivatives for 
 $\nu=1/3$ and $\nu = 2/3$ [9,12]. 
 For $\nu  = 1/4$, the solution was deduced by Barkai [13] in 2001 
  using Mathematica as a sum of three generalized hypergeometric functions, but the numerical result was uncertain, 
 presumably for a bag in the computing program . 
 In 2010-2012 Gorska and Penson [14]-[15] were able to represent 
 $f_{\nu,0}(t)$ in terms of Mejer $G$ functions.  
 Earlier, in 1958  Ragab [9]  expressed
 $f_{\nu,0}(t)$ in terms of   MacRobert  $E$ functions for  $\nu = 1/{n}$. 

 In 1952 Wlodarski [3] showed that when the generalized product Efros theorem [16] is applied to the Laplace transform of  
 $f_{\nu,\mu}(t)$  which is given in \eqref{GrindEQ__1_}, then it is possible to derive the following formula
\begin{equation} \label{GrindEQ__4_} 
\begin{array}{l} {\LL\left\{g(t)\right\}=G(s)\quad ;\quad 0<\nu <1\quad ;\quad \mu \ge 0} 
\\  \phantom{\rule{1pt}{15pt}}  
{\LL^{-1} \left\{\frac{G(s^{\nu } )}{s^{\mu } } \right\}={\ds \int _{0}^{\infty} }g(u)
 \LL^{-\, 1} \left\{F(s)\right\}  \, du=} 
\\  \phantom{\rule{1pt}{15pt}} 
{{\ds \int _{0}^{\infty }}g(u)  \LL^{-\, 1} \left\{\frac{\e^{-\, u s^{\nu } } }{s^{\mu } } \right\}  \, du={\ds \int _{0}^{\infty }}g(u) f_{\nu ,\mu } (u)  \, du}
 \end{array} 
\end{equation} 
It was established in our recent paper [17] that for
 $\mu = 0$ and $\mu = 1-\nu$ this functional expression can be written in terms
  of specific  Wright  functions $W_{\nu, \mu} (-t)$,  sometimes referred  to as  the Mainardi functions 
\begin{equation} \label{GrindEQ__5_} 
\begin{array}{l} {F_{\nu } (t)=W_{-\, \nu ,0} (-\, t)\quad ;\quad 0<\nu <1} 
\\  \phantom{\rule{1pt}{15pt}} 
{M_{\nu } (t)=W_{-\, \nu ,1\, -\, \nu } (-\, t)} \end{array} 
\end{equation} 
in the following way
\begin{equation} \label{GrindEQ__6_} 
\begin{array}{l} 
{\nu  \LL^{-1} \left\{\frac{G(s^{\nu } )}{s^{1\, -\, \nu } } \right\}=
{\ds \int _{0}^{\infty }}g(u)
 \, F_{\nu } \left(\frac{u}{t^{\nu } } \right)  \, \frac{du}{u} \quad ;\quad 0<\nu <1} 
\\  \phantom{\rule{1pt}{15pt}} 
 {t^{\nu }  \LL^{-1} \left\{\frac{G(s^{\nu } )}{s^{1\, -\, \nu } } \right\}
 = {\ds \int _{0}^{\infty }}g(u)
\,  M_{\nu } \left(\frac{u}{t^{\nu } } \right)  \, du} 
 \end{array} 
\end{equation} 
This follows from the fact that the functions
$F_\nu(t)$ and $M_\nu(t)$ 
satisfy
\begin{equation} \label{GrindEQ__7_} 
\begin{array}{l} {L\left\{\dfrac{1}{\lambda  \nu } F_{\nu } \left(\dfrac{\lambda }{t^{\nu } } \right)\right\}=L\left\{\dfrac{1}{t^{\nu } } 
M_{\nu } \left(\dfrac{\lambda }{t^{\nu } } \right)\right\}=\dfrac{e^{-\, \lambda \, s^{\nu } } }{s^{1\, -\, \nu } } }
 \\  \phantom{\rule{1pt}{15pt}} 
  {0<\nu <1\quad ;\quad \lambda >0} 
  \end{array} 
\end{equation} 
It was also illustrated by us in [17], that in many cases, by using standard tables of the Laplace transforms [18-21], the left hand side Laplace transforms in 
\eqref{GrindEQ__6_} can be inverted and numerous infinite integrals, 
finite integrals and integral identities for the 
 functions $F_\nu(t)$ and $M_\nu(t)$ can be derived.

\noindent 
Let us recall that the Wright functions [22-23], considered initially as a some kind generalization of the Bessel functions, are defined as 
an entire functions   of the argument  $z\in \CC$
 and parameters $\lambda>-1$  and $\mu \in \CC$  by
 \bee
W_{\lambda ,\mu } (z)={\ds \sum _{k = 0}^{\infty }\dfrac{z^{k} }{k!\, \, 
 \Gamma (\lambda   k+\mu )}} \,;\quad \lambda >- 1, \quad \mu \in \CC\,.
 \ee  
 It is usual to distinguish them in two kinds, 
 the first kind with $\lambda \ge0$ and the second kind with 
 $\lambda = -\nu$ and $\nu\in (0,1)$, see e.g [11].

Restricting  our attention to positive argument $ t>0$,
the Mainardi functions  turn out to be  particular Wright functions of the second kind  
 expressed by the following series
\begin{equation} \label{GrindEQ__9_} 
\begin{array}{l}
 F_{\nu } (t)=
{\ds  \sum _{k\, =\, 1}^{\infty }}\dfrac{(-\, t)^{k\, } }{k!\, \Gamma (-\, \nu \, k)}
 =
  \dfrac{1}{\pi }
  {\ds  \sum _{k\, =\, 1}^{\infty }}\dfrac{(-1)^{k\, +\, 1} t^{k} }{k!\, }  \,
   \Gamma (\nu \, k+1)\, \sin (\pi \, \nu \, k)\,,
   \\   \phantom{\rule{1pt}{15pt}} 
   M_{\nu } (t)=
   {\ds \sum _{k\, =\, 0}^{\infty }}\dfrac{(-\, t)^{k\, } }{k!\, \Gamma \, 
   \left(-\, \nu  (\, k+1)+1\right)}
    =\dfrac{1}{\pi } 
    {\ds \sum _{k\, =\, 1}^{\infty }}
    \dfrac{(-\, t)^{k\, -\, 1} }{(k-1)!\, }  \, \Gamma (\nu \, k)\, \sin (\pi \, \nu \, k)\,,
     \\   \phantom{\rule{1pt}{15pt}} 
     {F_{\nu } (t)=\nu  \, t \, M_{\nu } (t)}\,. 
     \end{array} 
\end{equation} 
\newpage
The interest in the Wright functions of the second kind  comes from the fact that they play an important role in solution of the linear partial  differential equations 
of fractional order which describe a wide spectrum of  phenomena including 
probability distributions,  anomalous diffusion and diffusive waves [24-33].

 In the section 2 of this paper, by using the complex inversion formula in \eqref{GrindEQ__1_} and the Bromwich contour, the integral representation of  
 $f_{\nu, \mu}(t)$ with $0<\nu<1$ and $0\le \mu <1$
   is derived and some basic properties of this inverse transform are established. 
 The cases $\mu\ge 1$ can be dealt as well, see Appendix for details.

  The next two sections 3 and 4  are devoted to evaluation of infinite, finite and convolution integrals by using the Efros theorem in the Wlodarski form. 
  The integrands of these integrals or integral identities include the elementary functions (power, exponential, logarithmic, trigonometric and hyperbolic functions) and the special functions (the error function, Mittag-Leffler functions and the Volterra functions). 
  The last section provides concluding remarks. 
  
  In derivations of these integrals direct and inverse Laplace transforms which are taken from tables of transforms [18-21] are always presented in  mathematical expressions. 
    All mathematical operations and manipulations with elementary and special functions, integrals and transforms are formal and their validity is assured by considering the restrictions usually imposed in the operational calculus.

\noindent 

\noindent
\section{Integral representations of the inverse Laplace \\ transform 
of $s^{-\mu}\, \exp (-s^\nu)$}  

Before infinite integrals in \eqref{GrindEQ__6_} will be evaluated, it is of interest to derive 
 the function $f_{\nu,\mu}(t)$ {by performing the complex integration from \eqref{GrindEQ__1_}. In investigated case, the branch point of the integrand exists and is located at the origin ${s} = 0$  and therefore the equivalent Bromwich contour is plotted in Figure \ref{f1}.
 The closed contour of integration ABCDEFA includes the line AB, the arcs BC and FA of a 
circle of radius $R \to \infty$\ with center at origin, the arc DE of a circle of radius 
$r \to 0$  with center at origin,  and two parallel lines CD and EF. 
The the cut along the negative axis  ensures that 
$F(s)$  is a single-valued function. 
However, according to the Cauchy lemma the integrals along the arcs BC  and FA 
vanish as  $R\to   \infty$.

%

\begin{figure}[h]
\centering
\includegraphics[scale=0.4]{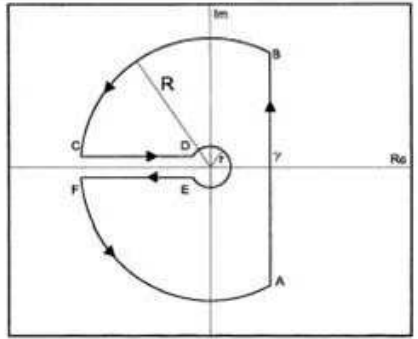}
\caption{The equivalent Bromwich contour.}
\label{f1}
\end{figure}

As a consequence 
there are only three contributions coming from integrals on the CD and EF 
lines and from the small circle round origin,D E, so we have
\begin{equation} \label{GrindEQ__10_} 
\begin{array}{l} {F(s)=\dfrac{\e^{- s^{\nu } } }{s^{\mu } } \quad ;\quad 0<\nu <1\quad ;\quad \mu \ge 0}, 
\\  \phantom{\rule{1pt}{15pt}}  
{f_{\nu ,\mu } (t)=
\dfrac{1}{2 \pi i} \, {\ds \int _{0}^{\infty }}\e^{-u t }
  \left[F(u \e^{-\, \pi  i} )-F(u \e^{\pi  i} )\right]\, du+
  \dfrac{1}{2 \pi i} {\mathop{\lim }\limits_{r \to 0 }} 
  {\ds \int _{-\, \pi }^{\pi }}\e^{t r \e^{i \theta } }  F(r \e^{i \theta } ) \, r \e^{i \theta } dr}.
 \end{array} 
\end{equation}

\noindent 
However, the last trigonometric integral vanishes for $\mu <1 $  and therefore the final result of the complex integration from \eqref{GrindEQ__10_} is
\begin{equation} \label{GrindEQ__11_} 
\begin{array}{l} 
{f_{\nu ,\mu } (t)=\dfrac{1}{\pi } \,
{\ds  \int _{0}^{\infty }}\dfrac{\e^{-u t\, -\, u^{\nu }  \cos (\pi  \nu ) } }{u^{\mu } }  
 \sin [  u^{\nu }  \sin (\pi  \nu )+\pi  \mu ]\, du}
 \\ \phantom{\rule{1pt}{15pt}}
  {t>0\quad ;\quad 0<\nu <1\quad ;\quad \mu <1} \end{array} 
\end{equation} 
The same integral representation has been derived in 1970 by Stankovi\.{c} [10] for the Wright function, but on the negative values of the argument $t$ and which is presented 
here  in our notation
\begin{equation} \label{GrindEQ__12_} 
\begin{array}{l} 
{t^{\mu \, -\, 1}  W_{\nu ,\mu } \left(-\dfrac{1}{t^{\nu } } \right)=
\dfrac{1}{\pi } \, {\ds \int _{0}^{\infty }}\dfrac{\e^{-u t\, -  u^{\nu }  \cos (\pi  \nu ) } }{u^{\mu } }
   \sin [  u^{\nu }  \sin (\pi  \nu )+\pi  \mu ]\, du} 
\\  \phantom{\rule{1pt}{15pt}}
{t>0\quad ;\quad 0<\nu <1\quad ;\quad \mu <1} \end{array} 
\end{equation} 
and comparing \eqref{GrindEQ__11_} with \eqref{GrindEQ__12_} we have that the inverse exponential functions can be expressed in terms of the Wright functions,
in agreement with the survey analysis by Mainardi and Consiglio [33] where also plots are presented, 
\begin{equation} \label{GrindEQ__13_} 
f_{\nu ,\mu } (t)=t^{\mu \, -\, 1}  W_{\nu ,\mu } \left(-\dfrac{1}{t^{\nu } } \right)\quad ;\quad 0<\nu <1\quad ;\quad \mu \ge 0 ,. 
\end{equation} 
In particular,  for $\mu  = 1 - \nu$ the Wright functions  are reduced to the Mainardi 
functions $F_\nu(t)$ and $M_\nu(t)$ [11] 
\begin{equation} \label{GrindEQ__14_} 
f_{\nu ,1\, -\, \nu } (t)
=
\dfrac{1}{ \nu } F_{\nu } \left(\dfrac{1}{t^{\nu } } \right)
=
\dfrac{1}{t^{\nu } } M_{\nu } \left(\dfrac{1}{t^{\nu } } \right)\quad ;\quad 0<\nu <1.
\end{equation} 
\newpage
For $\nu=1/2, \; \mu=0$ and 
$ \nu=\mu=1/2$, the integrals in \eqref{GrindEQ__11_} become the Laplace transforms of trigonometric functions [18-21]
\begin{equation} \label{GrindEQ__15_} 
\begin{array}{l} {f_{1/2,0} (t)
=
\dfrac{1}{\pi } \,{\ds  \int _{0}^{\infty }} \e^{-u t\,  }   \sin ( u^{1/2}  )\, du
=
\dfrac{\, \e^{-\, 1/4 t} }{2 \sqrt{\pi } \, t^{3/2} } }
 \\   \phantom{\rule{1pt}{15pt}}
 {f_{1/2,1/2} (t)=
 \dfrac{1}{\pi } \, {\ds \int _{0}^{\infty }}\dfrac{\e^{-u t\,  } }{u^{1/2} }
    \sin ( u^{1/2} +\dfrac{\pi }{2}  )\, du
    =\dfrac{1}{\pi } \, {\ds \int _{0}^{\infty }}\dfrac{\e^{-u t\,  } }{u^{1/2} }   \cos ( u^{1/2}  )\, du
    =\dfrac{\e^{-\, 1/4 t} }{\sqrt{\pi } } }
     \end{array} 
\end{equation} 
Since [11,12]
\begin{equation} \label{GrindEQ__16_} 
F_{1/3} \left(\dfrac{1}{t^{1/3} } \right)=
\dfrac{1}{3 \pi  t^{1/2} } \, K_{1/3} \left(\dfrac{2}{\sqrt{27 t} } \right) 
\end{equation} 
it follows from  \eqref{GrindEQ__11_} and \eqref{GrindEQ__14_} that
\begin{equation} \label{GrindEQ__17_} 
\begin{array}{l} {f_{1/3,2/3} (t)
=
\dfrac{1}{2 \pi } \, {\ds \int _{0}^{\infty }}\dfrac{e^{-u t\, -\, u^{1/3} /2 } }{u^{2/3} }   \left[\sqrt{3}  \cos \left(\dfrac{\sqrt{3}  }{2}   u^{1/3} \right)-  \sin \left(\dfrac{\sqrt{3}  }{2}  u^{1/3} \right)\right] \, du
=}
 \\ \phantom{\rule{1pt}{15pt}} 
 {\dfrac{1}{\pi } \sqrt{\dfrac{\lambda }{t} } \, K_{1/3} \left(\dfrac{2 }{\sqrt{27 t} } \right)} \end{array} 
\end{equation} 

 Differentiation of the integral \eqref{GrindEQ__11_} with respect to the argument \textit{t} gives
\begin{equation} \label{GrindEQ__18_} 
\begin{array}{l} 
{\dfrac{d ^{n} f_{\nu ,\mu } (t)}{d t^{n} } 
=
f_{\nu ,\mu \, -\, n} (t)=\dfrac{(-1)^{n} }{\pi } \,
{\ds  \int _{0}^{\infty }}\dfrac{\e^{-u t\, - u^{\nu }  \cos (\pi  \nu ) } }{u^{\mu \, -\, n} }   \sin [u^{\nu }  \sin (\pi  \nu )+\pi  \mu ]\, du} 
\\ \phantom{\rule{1pt}{15pt}}
 {t>0\quad ;\quad 0<\nu <1\quad ;\quad \mu <1\quad ;\quad n=0,1,2,3,...}
  \end{array} 
\end{equation} 
and using rules of the operational calculus we have the initial and final values of the function from [33]
\begin{equation} \label{GrindEQ__19_} 
\begin{array}{l} 
{f_{\nu ,\mu } (t\, \to \, +0)={\mathop{\lim }\limits_{s\, \to \, \infty }} [s F(s)]=s^{1\, -\, \mu } \e^{-\,  s^{\nu } } =0\quad ;\quad 0<\nu <1\quad ;\quad \mu <1}
 \\ \phantom{\rule{1pt}{15pt}}
  {f_{\nu ,\mu } (t\, \to \, \infty )={\mathop{\lim }\limits_{s\, \to \, 0}} [s F(s)]=s^{1\, -\, \mu } \e^{- s^{\nu } } =0} 
  \end{array} 
\end{equation} 
This is in an agreement with the finding of Pollard [2] that for $\mu = 0$, the inverse Laplace transform $f_{\nu,0}$ is positive almost everywhere, but 
for $\mu<0$  Stankovi\.{c} [10] postulated that in a some interval this function is negative and has at least one zero. 
Expanding the exponent in $F(s)$  into series it is possible to obtain the behaviour of the function for large values of the argument  ${t}$
\begin{equation} \label{GrindEQ__20_} 
\begin{array}{l} 
{ \LL^{-\, 1} \left\{F(s)\right\}=
\L^{-\, 1} \left\{\dfrac{e^{-\,  s^{\nu } } }{s^{\mu } } \right\}=
 \LL^{-\, 1} \left\{\dfrac{1}{s^{\mu } } \, \left[1-\dfrac{ s^{\nu  } }{1!} +\dfrac{ s^{2 \nu  } }
{2!} -...\right]\right\}} 
\\ \phantom{\rule{1pt}{15pt}} 
{f_{\nu ,\mu } (t\, \to \, \infty )\sim \dfrac{1}{t^{1\, -\, \mu } } \left\{\dfrac{1}{\Gamma (\mu ) } -\dfrac{1}{\Gamma (\mu -\nu )  t^{\nu } } +\dfrac{1}{2 \Gamma (\mu -2 \nu )  t^{2 \nu } } -...\right\}} 
\end{array} 
\end{equation} 
The integral of the 
$f_{\nu,\mu}$  can be derived from directly from
\begin{equation} \label{GrindEQ__21_} 
\int _{0}^{t}f_{\nu ,\mu } (u) \, du
=
\LL^{-\, 1} \left\{\dfrac{\e^{-\,  s^{\nu } } }{s^{1\, +\, \mu } } \right\}
=f_{\nu ,\mu \, +\, 1} (t) 
\end{equation} 
In order to obtain the recurrence relations the following operational rule can be applied
\begin{equation} \label{GrindEQ__22_} 
L\left\{t f_{\nu ,\mu \, -\, 1} (t)\right\}=-\dfrac{d}{d s} \left[F(s)\right]=-\dfrac{d}{d s} \left\{\dfrac{\e^{- s^{\nu } } }{s^{\mu \, -\, 1} } \right\}=(\mu -1)\dfrac{\e^{- s^{\nu } } }{s^{\mu } } +  \nu \dfrac{\e^{-\, s^{\nu } } }{s^{\mu \, -\, \nu } }  
\end{equation} 
and the inverse transforms in \eqref{GrindEQ__22_} are
\begin{equation} \label{GrindEQ__23_} 
t f_{\nu ,\mu \, -\, 1} (t)=(\mu -1) f_{\nu ,\mu } (t)+  \nu  f_{\nu ,\mu \, -\, \nu } (t).
\end{equation} 
However, from \eqref{GrindEQ__18_} with \textit{n} = 1, we have 
\begin{equation} \label{GrindEQ__24_} 
\dfrac{d f_{\nu ,\mu } (t)}{d t} =f_{\nu ,\mu \, -\, 1} (t) 
\end{equation} 
and therefore
\begin{equation} \label{GrindEQ__25_} 
t\, \dfrac{d f_{\nu ,\mu } (t)}{d t} =(\mu -1) f_{\nu ,\mu } (t)+  \nu  f_{\nu ,\mu \, -\, \nu } (t).
\end{equation}

\section{Integrals of the inverse Laplace transform  of \\ 
$s^{-\mu}\, \exp(s^\nu)$ with elementary functions}

\noindent In the first example the Wlodarski integral formula \eqref{GrindEQ__4_} is applied to the power function $g(t)=t^\lambda$,
\begin{equation} \label{GrindEQ__26_} 
\begin{array}{l} {g(t)=t^{\lambda } \quad ;\quad \lambda >0} 
\\ \phantom{\rule{1pt}{15pt}}  
{G(s)=\dfrac{\Gamma (\lambda +1)}{s^{\lambda \, +\, 1} } \, \quad ;\quad G(s^{\nu } )=\dfrac{\Gamma (\lambda +1)}{s^{(\lambda \, +\, 1)  \nu } } } 
\\ \phantom{\rule{1pt}{15pt}}  
{\LL^{-1} \left\{\dfrac{1}{s^{\mu } } \bullet \dfrac{\Gamma (\lambda +1)}{s^{(\lambda \, +\, 1)  \nu } } \right\}
= \LL^{-1} \, \left\{\dfrac{\Gamma (\lambda +1)}{s^{(\lambda \, +\, 1)  \nu \, +\, \mu } } \right\}=\dfrac{\Gamma (\lambda +1) \, t^{(\lambda \, +\, 1)  \nu \, +\mu \, -\, 1} }{\Gamma   [ (\lambda \, +\, 1)  \nu \, +\mu )]} } \end{array} 
\end{equation} 
and therefore we have
\begin{equation} \label{GrindEQ__27_} 
{\ds \int _{0}^{\infty }} u^{\lambda \, } f_{\nu ,\mu } (u) \,
  du=\dfrac{\Gamma (\lambda +1) \, t^{(\lambda \, +\, 1)  \nu \, +\mu \, -\, 1} }{\Gamma   [ (\lambda \, +\, 1)  \nu \, +\mu )]} \quad ;\, 0<\nu <1\quad ;\, \mu <1\quad ;\, \lambda >0. 
\end{equation}

 The above results can be extended to functions $g(t)$ which are defined at finite intervals and to step functions jumping at integral values of  variable $t$. Let start with 
\begin{equation} \label{GrindEQ__28_} 
\begin{array}{l} {g(t)=\left\{\begin{array}{l} {1\quad ;\quad 0<t<\lambda }
 \\ {0\quad ;\quad t>\lambda } \end{array}\right. }
  \\ \phantom{\rule{1pt}{15pt}}   
  {G(s)=\dfrac{1-\e^{-\, \lambda  s} }{s} \, \quad ;\quad G(s^{\nu } )=
  \dfrac{1-\e^{-\, \lambda  s^{\nu } } }{s^{  \nu } } }
 \\ \phantom{\rule{1pt}{15pt}}  
{\LL^{-1} \left\{\dfrac{1}{s^{\mu } } \bullet \dfrac{1-\e^{-\, \lambda  s^{\nu } } }{s^{  \nu } } \right\}
= 
\LL^{-1} \, \left\{\dfrac{1}{s^{\nu \, +\, \mu } } -\dfrac{\lambda ^{(\mu /\nu \, +\, 1)}
  \e^{-\, (\lambda ^{1/\nu }  s)^{\nu } } }{(\lambda ^{1/\nu } s)^{\nu \, +\, \mu } } \right\}=} 
\\  \phantom{\rule{1pt}{15pt}}  
{\dfrac{t^{\nu \, +\, \mu \, -\, 1} }{\Gamma (\nu \, +\, \mu )} -\lambda ^{[  (\mu \; -\, 1)/\nu \, +\, 1]}  f_{\nu ,\nu \, +\, \mu } \left(\dfrac{t}{\lambda ^{1/\nu } } \right)} 
\end{array} 
\end{equation} 
which leads to the finite integral
\begin{equation} \label{GrindEQ__29_} 
\begin{array}{l} {
{\ds \int _{0}^{\lambda }}f_{\nu ,\mu } (u)  \, du
=
\dfrac{t^{\nu \, +\, \mu \, -\, 1} }
{\Gamma (\nu \, +\, \mu )} -\lambda ^{[  (\mu \; -\, 1)/\nu \, +\, 1]}  f_{\nu ,\nu \, +\, \mu } (\frac{t}{\lambda ^{1/\nu } } )} 
\\  \phantom{\rule{1pt}{15pt}}  
{0<\nu <1\quad ;\, \mu <1\quad ;\, \lambda >0} 
\end{array} 
\end{equation} 
In the next example, from
\begin{equation} \label{GrindEQ__30_} 
\begin{array}{l} {g(t)=
\left\{\begin{array}{l} 
{1-t\quad ;\quad 0<t<1} \\ {0\quad ;\quad t>1} \end{array}\right. \quad } 
\\ \phantom{\rule{1pt}{15pt}}  
 {G(s)=\dfrac{\e^{- s} +s-1}{s^{2} } \, \quad ;\quad G(s^{\nu } 
 )=\dfrac{\e^{- s^{\nu } } +s^{\nu } -1}{s^{  2 \nu } } }
 \\ \phantom{\rule{1pt}{15pt}}  
  {\LL^{-1} \left\{\dfrac{1}{s^{\mu } } \bullet \dfrac{\e^{- s^{\nu } } +s^{\nu } -1}{s^{2 \nu } } \right\}
  = \LL^{-1} \, \left\{\dfrac{1}{s^{2 \nu \, +\, \mu } } -\dfrac{1}{s^{\nu \, +\, \mu } } +\dfrac{\e^{-\,  s^{\nu } } }{s^{2 \nu \, +\, \mu } } \right\}=} 
  \\ 
  {\dfrac{t^{2 \nu \, +\, \mu \, -\, 1} }{\Gamma (2 \nu \, +\, \mu )} -\dfrac{t^{ \nu \, +\, \mu \, -\, 1} }{\Gamma (\nu \, +\, \mu )} +f_{\nu ,2 \nu \, +\, \mu } (t)} 
  \end{array} 
\end{equation} 
we have 
\begin{equation} \label{GrindEQ__31_} 
{\ds \int _{0}^{1}}(1-u) f_{\nu ,\mu } (u)  \, du
=\dfrac{t^{2 \nu \, +\, \mu \, -\, 1} }{\Gamma (2 \nu \, +\, \mu )} -\dfrac{t^{ \nu \, +\, \mu \, -\, 1} }{\Gamma (\nu \, +\, \mu )} +f_{\nu ,2 \nu \, +\, \mu } (t) 
\end{equation} 
Using
\begin{equation} \label{GrindEQ__32_} 
\begin{array}{l}
 {g(t)=\left\{\begin{array}{l} {0\quad ;\quad 0<t<\lambda } \\ {\dfrac{(t-\lambda )^{\mu } }{\Gamma (\mu +1)} \quad ;\quad t>\lambda } \end{array}\right. \quad } 
 \\ \phantom{\rule{1pt}{15pt}}  
  {G(s)=\dfrac{\e^{- \lambda s} }{s^{\mu \, +\, 1} } \, \quad ;\quad G(s^{\nu } )=\dfrac{e^{- \lambda  s^{\nu } } }{s^{(\mu \, +\, 1) \nu } } }
 \\ \phantom{\rule{1pt}{15pt}}  
  {\LL^{-1} \left\{\dfrac{1}{s^{\mu } } \bullet \dfrac{\e^{- \lambda  s^{\nu } } }{s^{(\mu \, +\, 1) \nu } } \right\}=
  \LL^{-1} \left\{\dfrac{\lambda ^{[  (1/\nu \, +\, 1) \mu \, +\, 1]}  \e^{-\, (\lambda ^{1/\nu }  s)^{\nu } } }{(\lambda ^{1/\nu } s)^{(\mu \, +\, 1) \nu \, +\, \mu } } \right\}=} 
  \\ \phantom{\rule{1pt}{15pt}}  
  {=\lambda ^{[  (1/\nu \, +\, 1) \mu \, \, -\, 1/\nu \, +\, 1]} f_{\nu ,\nu \, +\, \nu  \mu \, + \mu } (\frac{t}{\lambda ^{1/\nu } } )}
   \end{array} 
\end{equation} 
it is possible to derive
\begin{equation} \label{GrindEQ__33_} 
{\ds \int _{\lambda }^{\infty} }(t-\lambda )^{\mu }  f_{\nu ,\mu } (u)   \, du
=
\lambda ^{[(1/\nu \, +\, 1) \mu \, \, -\, 1/\nu \, +\, 1]} f_{\nu ,\nu \, +\, \nu  \mu \, + \mu } \left(\frac{t}{\lambda ^{1/\nu } } \right) .
\end{equation} 
 In the next two examples of the case of exponential functions are considered
\begin{equation} \label{GrindEQ__34_} 
\begin{array}{l} 
{g(t)=\e^{-\, \alpha   t} \quad ;\quad \alpha >\, 0} 
\\  \phantom{\rule{1pt}{15pt}}  
{G(s)=\dfrac{1}{s+\alpha } \quad ;\, \quad G(s^{\nu } )=\dfrac{1}{s^{\nu } +\alpha } } 
\\ \phantom{\rule{1pt}{15pt}}  
{\LL^{-1} \left\{\dfrac{1}{s^{\mu } } \bullet \dfrac{1}{s^{\nu } +\alpha } \right\}
=
\LL^{-1} \, \left\{\dfrac{1}{s^{\nu \, +\mu \, -1} } \bullet \dfrac{s^{\nu \, -\, 1} }{s^{\nu } +\alpha } \right\}=
\dfrac{t^{\nu \, +\, \mu \, -2} }{\Gamma (\nu \, +\, \mu \, -1)}
\, \star  \, E_{\nu } (-\, \alpha   t^{\nu } )} 
\end{array} 
\end{equation} 
From \eqref{GrindEQ__4_} and \eqref{GrindEQ__34_} the convolution integral with the Mittag-Leffler function
\begin{equation} \label{GrindEQ__35_} 
\begin{array}{l}
 {E_{\nu } (t)={\ds \sum _{k\, =\, 0}^{\infty }}\dfrac{t^{k} }{\Gamma (k \nu +1)}\,,  } 
\\ \phantom{\rule{1pt}{15pt}}  
{\LL\left\{E_{\nu } (-\alpha   t^{\nu } )\right\}=\dfrac{s^{\nu \, -\, 1} }{s^{\nu } +\alpha }\,, } \end{array} 
\end{equation} 
is derived
\begin{equation} \label{GrindEQ__36_} 
\begin{array}{l} 
{{\ds \int _{0}^{\infty }}\e^{\, -\, } {}^{\alpha  u\, } f_{\nu ,\mu } (u)   du
=\dfrac{t^{\nu \, +\, \mu \, -2} }{\Gamma (\nu \, +\, \mu \, -1)} 
\, \star \, E_{\nu } (-\, \alpha   t^{\nu } )}
 \\ \phantom{\rule{1pt}{15pt}}  
  {0<\nu <1\quad ;\quad \mu <1\quad ;\quad \alpha >0\,.} 
  \end{array} 
\end{equation} 
By changing variable of integration $x = t\,  (cos\theta)^2$, all convolution integrals can be expressed as in terms of finite trigonometric integrals. The shifted increasing and decreasing exponential functions are considered in the following two examples. 
From
\begin{equation} \label{GrindEQ__37_} 
\begin{array}{l}
 {g(t)=\left\{\begin{array}{l} {0\quad ;\quad 0<t<\lambda \quad ;\quad 0<\nu <1} 
\\ 
{1-\e^{-\, (t\, -\, \lambda )} \quad ;\quad t>\lambda } \end{array}\right. \quad } 
\\ \phantom{\rule{1pt}{15pt}}  
{G(s)=\dfrac{e^{- \lambda s} }{s (s+1)} \, \quad ;\quad G(s^{\nu } )
=
\dfrac{\e^{- \lambda  s^{\nu } } }{s^{ \nu } (s^{ \nu } +1)} } 
\\ \phantom{\rule{1pt}{15pt}}  
{\LL^{-1} \left\{\dfrac{1}{s^{\mu } } \bullet \dfrac{\e^{- \lambda  s^{\nu } } 
}{s^{ \nu } (s^{ \nu } +1)} \right\}
=
\LL^{-1} \, \left\{\dfrac{s^{ \nu \, -\, 1} }{(s^{ \nu } +1)} \bullet \dfrac{\lambda ^{(2 \nu \, +\, \mu \, -\, 1)/\nu } \e^{-\, (\lambda ^{1/\nu }  s)^{\nu } } }{(\lambda  ^{1/\nu } s)^{(2 \nu \, +\, \mu \, -\, 1)} } \right\}=} 
\\ \phantom{\rule{1pt}{15pt}}  
 {\lambda ^{(\nu \, +\, \mu \, -\, 1)/\nu } E_{\nu } (- \, t^{\nu } )*f_{\nu ,} {}_{2 \nu \, +\, \mu \, -\, 1} (\frac{t}{\lambda  ^{1/\nu } } )}
  \end{array} 
\end{equation} 
we have the convolution of two functions.
\begin{equation} \label{GrindEQ__38_} 
{\ds \int _{\lambda }^{\infty }}
\left [1-\e^{-\, (u\, -\, \lambda )} \right]  f_{\nu ,\mu } (u)  \, du
=
\lambda ^{(\nu \, +\, \mu \, -\, 1)/\nu } E_{\nu } (- \, t^{\nu } )
\, \star \, f_{\nu ,} {}_{2 \nu \, +\, \mu \, -\, 1} (\frac{t}{\lambda  ^{1/\nu } } ). 
\end{equation} 
Similarly from
\begin{equation} \label{GrindEQ__39_} 
\begin{array}{l} 
{g(t)=\left\{\begin{array}{l} {0\quad ;\quad 0<t<\lambda \quad ;\quad 0<\nu <1}
\\ {\e^{-\, (t\, -\, \lambda )} \quad ;\quad t>\lambda } \end{array}\right. \quad }
 \\ \phantom{\rule{1pt}{15pt}}  
  {G(s)=\dfrac{e^{- \lambda s} }{ (s+1)} \, \quad ;\quad G(s^{\nu } )=\dfrac{e^{- \lambda  s^{\nu } } }{(s^{ \nu } +1)} } 
  \\ \phantom{\rule{1pt}{15pt}}   
{\LL^{-1} \left\{\dfrac{1}{s^{\mu } } \bullet \dfrac{\e^{- \lambda  s^{\nu } } }{(s^{ \nu } +1)} \right\}=
\LL^{-1} \, \left\{\dfrac{s^{ \nu \, -\, 1} }{(s^{ \nu } +1)} \bullet \dfrac{\lambda ^{(\nu \, +\, \mu \, -\, 1)/\nu } \e^{-\, (\lambda ^{1/\nu }  s)^{\nu } } }{(\lambda ^{1/\nu } s)^{(\nu \, +\, \mu \, -\, 1)}  } \right\}=} 
\\ \phantom{\rule{1pt}{15pt}}  
 {\lambda ^{(\nu \, +\, \mu \, -\, 2)/\nu } E_{\nu } (- \, t^{\nu } )*f_{\nu ,} {}_{\nu \, +\, \mu \, -\, 1} (\frac{t}{\lambda ^{1/\nu } } )} 
 \end{array} 
\end{equation} 
it is possible to obtain
\begin{equation} \label{GrindEQ__40_} 
{\ds \int _{\lambda }^{\infty }}\e^{-\, (u\, -\, \lambda )} \, f_{\nu ,\mu } (u)  \, du
=\lambda ^{(\nu \, +\, \mu \, -\, 2)/\nu } \, E_{\nu } (- \, t^{\nu } )
\, \star \, f_{\nu ,} {}_{\nu \, +\, \mu \, -\, 1} (\frac{t}{\lambda ^{1/\nu } } ). 
\end{equation} 

 The logarithmic functions are the next group of elementary functions to be considered. In the simplest case from
\begin{equation} \label{GrindEQ__41_} 
\begin{array}{l} 
{g(t)=\ln t\quad ;\quad C=\e^{\gamma } }
 \\ \phantom{\rule{1pt}{15pt}}  
 {G(s)=-\dfrac{\ln (C s)}{s} \, \quad ;\quad G(s^{\nu } )=-\dfrac{\ln (C s^{\nu } )}{s^{\nu } } } \\ \phantom{\rule{1pt}{15pt}}  
  { \LL^{-1} \left\{\dfrac{1}{s^{\mu } } \bullet \dfrac{-\ln (C s^{\nu } )}{s^{\nu } } \right\}= 
  \LL^{-1} \, \left\{\dfrac{1}{s^{\nu \, +\, \mu \, -\, 1} } \bullet \dfrac{-  \nu \ln (C s)+(\nu -1) \gamma }{s} \right\}=} 
  \\ \phantom{\rule{1pt}{15pt}}  
  {\dfrac{(\nu -1)  \gamma  t^{\nu \, +\, \mu \, -\, 1} }{\Gamma (\nu \, +\, \mu \, )} +\dfrac{\nu   t^{\nu \, +\, \mu \, -\, 2} }{\Gamma (\nu \, +\, \mu \, -\, 1)}   *\ln t} \end{array}           
\end{equation} 
where $\gamma$  is the Euler constant,
it follows that
\begin{equation} \label{GrindEQ__42_} 
{\ds \int _{0}^{\infty} }\ln u \, f_{\nu ,\mu } (u)  \, du
=
\dfrac{(\nu -1)  \gamma  t^{\nu \, +\, \mu \, -\, 1} }
{\Gamma (\nu \, +\, \mu \, )} +\dfrac{\nu   t^{\nu \, +\, \mu \, -\, 2} }{\Gamma (\nu \, +\, \mu \, -\, 1)}   \, \star \, \ln t\,.
\end{equation}

\noindent In the more general case
\begin{equation} \label{GrindEQ__43_} 
\begin{array}{l} 
{g(t)=t^{\lambda \, -\, 1} \ln t\quad ;\quad \lambda >0\quad ;\quad 0<\nu <1} 
\\ \phantom{\rule{1pt}{15pt}}  
{G(s)=\dfrac{\Gamma (\lambda )}{s^{\lambda } } \left[\psi (\lambda )-\ln s\, \right]\quad ;\quad \quad G(s^{\nu } )=\dfrac{\Gamma (\lambda )\, [\psi (\lambda )-\nu \ln s]}{s^{\lambda  \nu } } }
 \\ \phantom{\rule{1pt}{15pt}}  
 { \Gamma (\lambda )\, \LL^{-1} \left\{\dfrac{1}{s^{\mu } } \bullet \dfrac{\, [\psi (\lambda )-\nu \ln s]}{s^{\lambda  \nu } } \right\}=
 \Gamma (\lambda )\, 
  \LL^{-1} \, \left\{\dfrac{\psi (\lambda )+\gamma  \nu }{s^{(\lambda \,  \nu \, +\, \mu )} } -\dfrac{\nu }{s^{(\lambda \,  \nu \, +\, \mu \, -\, 1)} } \bullet \dfrac{\ln (C s)}{s} \right\}=} 
  \\ \phantom{\rule{1pt}{15pt}}   
  {\dfrac{\Gamma (\lambda ) [\psi (\lambda )+\gamma  \nu ]\, t^{(\lambda \,  \nu \, +\, \mu \, -1)} \, }{\Gamma (\lambda \,  \nu \, +\, \mu )} \, +\dfrac{ \Gamma (\lambda )\, \nu   t^{(\lambda \,  \nu \, +\, \mu \, -2)\, } }{\Gamma (\lambda \,  \nu \, +\, \mu \, -1)} \, \star \,\ln t^{\nu } } \end{array} 
\end{equation} 
we have
\begin{equation} \label{GrindEQ__44_} 
{\ds \int _{0}^{\infty }}u^{\lambda \, -\, 1} \ln u \, f_{\nu ,\mu } (u)  \, du
=\dfrac{\Gamma (\lambda )
\left [\psi (\lambda )+\gamma  \nu \right]
\, t^{(\lambda \,  \nu \, +\, \mu \, -1)} \, }{\Gamma (\lambda \,  \nu \, +\, \mu )} \,
 +\dfrac{ \Gamma (\lambda )\, \nu   t^{(\lambda \,  \nu \, +\, \mu \, -2)\, } }{\Gamma (\lambda \,  \nu \, +\, \mu \, -1)} \, \star \, \ln t^{\nu }  
\end{equation} 

 Trigonometric and hyperbolic functions is the last groups of elementary functions to be considered. From
\begin{equation} \label{GrindEQ__45_} 
\begin{array}{l} 
{g(t)=\sin (\lambda   t)\quad ;\quad 0<\nu <1\quad } 
\\ \phantom{\rule{1pt}{15pt}}  
{G(s)=\dfrac{\lambda }{s^{2} +\lambda ^{2} } \, \quad ;\quad G(s^{\nu } ) 
=
\dfrac{\lambda }{s^{2 \nu } +\lambda ^{2} } \, } 
\\ \phantom{\rule{1pt}{15pt}}  
{\lambda  \LL^{-1} \left\{\dfrac{1}{s^{\mu } } \bullet \dfrac{1}{s^{2 \nu } +\lambda ^{2} } \right\}=
\lambda  
\LL^{-1} \, \left\{\dfrac{1}{s^{2 \nu \, +\, \mu \, -\, 1} } 
\bullet \dfrac{s^{2 \nu \, -\, 1} }{s^{2 \nu } +\lambda ^{2} } \right\}
=
\dfrac{\lambda  t^{2 \nu \, +\, \mu \, -\, \, 2} }{\Gamma (2 \nu \, +\, \mu \, -\, 1)} 
\, \star \, E_{2 \nu } (-\lambda ^{2} \, t^{2 \nu } )},
 \end{array} 
\end{equation} 
the following integral identity is derived 
\begin{equation} \label{GrindEQ__46_} 
{\ds \int _{0}^{\infty} }\sin (\lambda  u)\, f_{\nu ,\mu } (u)  \, du=\dfrac{\lambda  t^{2 \nu \, +\, \mu \, -\, \, 2} }{\Gamma (2 \nu \, +\, \mu \, -\, 1)} \, \star\,
E_{2 \nu } (-\lambda ^{2} t^{2 \nu } ) .
\end{equation} 
\newpage
Similarly as in \eqref{GrindEQ__45_}, for the hyperbolic sine function,  the change is only in the sign
\begin{equation} \label{GrindEQ__47_} 
\begin{array}{l}
 {g(t)=\sinh (\lambda   t)\quad ;\quad 0<\nu <1\quad }
  \\ \phantom{\rule{1pt}{15pt}}  
   {G(s)=\dfrac{\lambda }{s^{2} -\lambda ^{2} } \, \quad ;\quad G(s^{\nu } )
   =
   \dfrac{\lambda }{s^{2 \nu } -\lambda ^{2} } \, } 
   \\ \phantom{\rule{1pt}{15pt}}  
   {\lambda  \LL^{-1} \left\{\dfrac{1}{s^{\mu } } \bullet \dfrac{1}{s^{2 \nu } -\lambda ^{2} } \right\}
   =
   \lambda  \LL^{-1} \, \left\{\dfrac{1}{s^{2 \nu \, +\, \mu \, -\, 1} } \bullet \dfrac{s^{2 \nu \, -\, 1} }{s^{2 \nu } -\lambda ^{2} } \right\}=\dfrac{\lambda  t^{2 \nu \, +\, \mu \, -\, \, 2} }{\Gamma (2 \nu \, +\, \mu \, -\, 1)} \, \star\, E_{2 \nu } (\lambda ^{2} t^{2 \nu } )}
    \end{array} 
\end{equation} 
and therefore we have
\begin{equation} \label{GrindEQ__48_} 
{\ds \int _{0}^{\infty }}\sinh (\lambda  u)\, f_{\nu ,\mu } (u)  \, du=
 \dfrac{\lambda  t^{2 \nu \, +\, \mu \, -\, \, 2} }{\Gamma (2 \nu \, +\, \mu \, -\, 1)} *E_{2 \nu } (\lambda ^{2} t^{2 \nu } ) .
\end{equation} 
In the case of cosine function, from 
\begin{equation} \label{GrindEQ__49_} 
\begin{array}{l} 
{g(t)=\cos (\lambda   t)\quad ;\quad 0<\nu <1}
 \\  \phantom{\rule{1pt}{15pt}}  
  {G(s)=\dfrac{s}{s^{2} +\lambda ^{2} } \, \quad ;\quad G(s^{\nu } )=\dfrac{s^{\nu } }{s^{2 \nu } +\lambda ^{2} } }
  \\ \phantom{\rule{1pt}{15pt}}  
   {\LL^{-1} \left\{\dfrac{1}{s^{\mu } } \bullet \dfrac{s^{\nu } }{s^{2 \nu } +\lambda ^{2} } \right\}
   = \LL^{-1} \, \left\{\dfrac{1}{s^{\nu \, +\, \mu \, -\, 1} } \bullet \dfrac{s^{2 \nu \, -\, 1} }{s^{2 \nu } +\lambda ^{2} } \right\}=
   \dfrac{t^{\nu \, +\, \mu \, -\, 2} }{\Gamma (\nu \, +\, \mu \, -1)}
   \, \star \, E_{2 \nu } (-\lambda ^{2} t^{2 \nu } ),} 
   \end{array} 
\end{equation} 
we have
\begin{equation} \label{GrindEQ__50_} 
{\ds \int _{0}^{\infty} }\cos (\lambda  u)\, f_{\nu ,\mu } (u)  \, du=\dfrac{t^{\nu \, +\, \mu \, -\, 2} }{\Gamma (\nu \, +\, \mu \, -1)} \,\star\, E_{2 \nu } (-\lambda ^{2} t^{2 \nu } ). 
\end{equation} 
Similarly as in \eqref{GrindEQ__47_} and \eqref{GrindEQ__48_}
\begin{equation} \label{GrindEQ__51_} 
\begin{array}{l} 
{g(t)=\cosh (\lambda   t)\quad ;\quad 0<\nu <1} 
\\ \phantom{\rule{1pt}{15pt}}  
{G(s)=\dfrac{s}{s^{2} -\lambda ^{2} } \, \quad ;\quad G(s^{\nu } )
=
\dfrac{s^{\nu } }{s^{2 \nu } -\lambda ^{2} } } 
\\ \phantom{\rule{1pt}{15pt}}  
{ \LL^{-1} \left\{\dfrac{1}{s^{\mu } } \bullet \dfrac{s^{\nu } }{s^{2 \nu } -\lambda ^{2} } \right\}= 
\LL^{-1} \, \left\{\dfrac{1}{s^{\nu \, +\, \mu \, -\, 1} } \bullet \dfrac{s^{2 \nu \, -\, 1} }{s^{2 \nu } -\lambda ^{2} } \right\}=\dfrac{t^{\nu \, +\, \mu \, -\, 2} }{\Gamma (\nu \, +\, \mu \, -1)} \, \star \, E_{2 \nu } (\lambda ^{2} t^{2 \nu } ),} 
\end{array} 
\end{equation} 
for the hyperbolic cosine function we have
\begin{equation} \label{GrindEQ__52_} 
{\ds \int _{0}^{\infty} }\cosh (\lambda  u)\, f_{\nu ,\mu } (u)  \, du
=
\dfrac{t^{\nu \, +\, \mu \, -\, 2} }{\Gamma (\nu \, +\, \mu \, -1)} 
\, \star \, E_{2 \nu } (\lambda ^{2} t^{2 \nu } ). 
\end{equation} 
In the case of the product of trigonometric and hyperbolic sine functions the direct and inverse Laplace transforms are  
\begin{equation} \label{GrindEQ__53_} 
\begin{array}{l} 
{g(t)=\sin (\lambda  t)\sinh (\lambda  t)\quad ;\quad 0<\nu <1} 
\\ \phantom{\rule{1pt}{15pt}}  
 {G(s)=\dfrac{2 \lambda ^{2} s}{s^{4} +4\lambda ^{4} } \, \quad ;\quad G(s^{\nu } )
 =
 \dfrac{2 \lambda ^{2} s^{\nu } }{s^{4 \nu } +4\lambda ^{4} } \, } 
 \\ \phantom{\rule{1pt}{15pt}}  
 { 2 \lambda ^{2} \LL^{-1} \left\{\dfrac{1}{s^{\mu } } \bullet \dfrac{s^{\nu } }{s^{4 \nu } +4\lambda ^{4} } \right\}
 = 2 \lambda ^{2} \LL^{-1} \, \left\{\dfrac{1}{s^{3 \nu \, +\, \mu \, -\, 1} } \bullet \dfrac{s^{4  \nu \, -\, 1} }{(s^{4 \nu } +4 \lambda ^{4} )} \right\}=}
  \\   \phantom{\rule{1pt}{15pt}}  
  {\frac{2 \lambda ^{2} }{\Gamma (3 \nu \, +\, \mu \, -1)}    [t^{(3 \nu \, +\, \mu \, -2)} 
  \, \star \, E_{4 \nu } (-\, 4\lambda ^{4} t^{4 \nu } )]} 
  \end{array} 
\end{equation} 
and therefore
\begin{equation} \label{GrindEQ__54_} 
\int _{0}^{\infty }\sin (\lambda  u)  \sinh (\lambda  u)\,f_{\nu ,\mu } (u)  \, du
=
\frac{2 \lambda ^{2} }{\Gamma (3 \nu \, +\, \mu \, -1)}    [t^{(3 \nu \, +\, \mu \, -2)}
\, \star \, E_{4 \nu } (-\, 4\lambda ^{4} t^{4 \nu } )]. 
\end{equation} 
For the product of trigonometric and hyperbolic sine functions we have 
\begin{equation} \label{GrindEQ__55_} 
\begin{array}{l} 
{g(t)=\cos (\lambda  t)\cosh (\lambda  t)\quad ;\quad 0<\nu <1}
 \\ \phantom{\rule{1pt}{15pt}}   
 {G(s)=\dfrac{s^{3 } }{s^{4} +4\lambda ^{4} } \, \quad ;\quad G(s^{\nu } )
 =\dfrac{s^{3 \nu } }{s^{4 \nu } +4\lambda ^{4} } \, } 
 \\  \phantom{\rule{1pt}{15pt}}  
 { \LL^{-1} \left\{\dfrac{1}{s^{\mu } } \bullet \dfrac{s^{3 \nu } }{s^{4} +4\lambda ^{4} } \right\}=
  \LL^{-1} \, \left\{\dfrac{1}{s^{\nu \, +\, \mu \, -\, 1} } \bullet \dfrac{s^{4  \nu \, -\, 1} }{s^{4 \nu } +4 \lambda ^{4} } \right\}
  =
  \dfrac{t^{\nu \, +\, \mu \, -\, 2} }{\Gamma (\nu \, +\, \mu \, -1)}
  \, \star \, E_{4 \nu } (-\, 4\lambda ^{4} t^{4 \nu } ),}
   \end{array} 
\end{equation} 
and therefore 
\begin{equation} \label{GrindEQ__56_} 
{\ds \int _{0}^{\infty} }\cos (\lambda  u) (\cosh (\lambda  u)\, f_{\nu ,\mu } (u)  \, du
=
\frac{t^{\nu \, +\, \mu \, -\, 2} }{\Gamma (\nu \, +\, \mu \, -1)}
\, \star  \, E_{4 \nu } (-\, 4\lambda ^{4} t^{4 \nu } ). 
\end{equation} 
\textbf{}

\section{Integrals of the inverse Laplace transform  of \\
$s^{-\mu}\, \exp(-s^\nu)$  with the Mitag-Leffler,  Error  and  \\ Volterra functions.}

\noindent The Laplace transform of the two parameter Mittag-Leffler function is  [11,34,35]
\begin{equation} \label{GrindEQ__57_} 
\begin{array}{l} 
{\LL\left\{t^{\beta \, -\, 1}\, E_{\alpha ,\beta } (\pm   \lambda   t^{\alpha } )\right\}=\dfrac{s^{\alpha \, -\, \beta } }{s^{\alpha } \mp \lambda }, } 
\\ \phantom{\rule{1pt}{15pt}}  
{E_{\alpha ,\beta } (z)=
{\ds \sum _{k\, =\, 0}^{\infty }}\dfrac{z^{k} }{\Gamma (k \alpha +\beta )},  } 
\end{array} 
\end{equation} 
and this permits to obtain 
\begin{equation} \label{GrindEQ__58_} 
\begin{array}{l} 
{g(t)=t^{\beta \, -\, 1} E_{\alpha ,\beta } (\pm   \lambda   t^{\alpha } )\quad ;\quad 0<\nu <1} 
\\ \phantom{\rule{1pt}{15pt}}  
 {G(s)=\dfrac{s^{\alpha \, -\, \beta } }{s^{\alpha } \mp \lambda } \quad ;\quad 
 G(s^{\nu } )=\dfrac{s^{(\alpha \, -\, \beta ) \nu } }{s^{\alpha  \nu } \mp \lambda } \quad }
  \\ \phantom{\rule{1pt}{15pt}}  
   {\left\{\dfrac{1}{s^{\mu } } \bullet \dfrac{s^{(\alpha \, -\, \beta ) \nu } }{s^{\alpha  \nu } \mp \lambda } \right\}=
    \LL^{-1} \left\{\dfrac{s^{\alpha  \nu \, -\, (\beta \,  \nu \, +\mu )} }{s^{\alpha  \nu } \mp \lambda } \right\}
    =
    t^{(\beta \,  \nu \, +\mu \, -\, 1) } E_{\alpha  \nu ,\beta \,  \nu \, +\, \mu } (\pm   \lambda   t^{\alpha } ),}
     \end{array} 
\end{equation} 
and
\begin{equation} \label{GrindEQ__59_} 
{\ds \int _{0}^{\infty }}u^{\beta \, -\, 1}
 E_{\alpha ,\beta } (\pm   \lambda   u^{\alpha } )  
\,f_{\nu ,\mu } (u)  \, du=\, t^{(\beta \,  \nu \, +\mu \, -\, 1) } E_{\alpha  \nu ,\beta \,  \nu \, +\, \mu } (\pm   \lambda   t^{\alpha } ). 
\end{equation} 
Thus, in both sides of expressions \eqref{GrindEQ__59_} appear the Mittag-Leffler functions. Evidently, for  $\beta= 1$  they are reduced to the classical Mittag-Leffler functions. 
\begin{equation} \label{GrindEQ__60_} 
{\ds \int _{0}^{\infty} }E_{\alpha ,1} (\pm   \lambda   u^{\alpha } )  \, f_{\nu ,\mu } (u)  \, du
= t^{(\,  \nu \, +\mu \, -\, 1) }\, E_{\alpha  \nu ,\,  \nu \, +\, \mu } (\pm   \lambda   t^{\alpha } ), 
\end{equation} 
 and for $\beta= \alpha$ and $\beta = \alpha  + 1$ we have
\begin{equation} \label{GrindEQ__61_} 
\begin{array}{l} 
{{\ds \int _{0}^{\infty }}u^{\alpha \, -\, 1} E_{\alpha ,\alpha } (\pm   \lambda   u^{\alpha } )  \, f_{\nu ,\mu } (u)  \, du =
 t^{(\alpha \,  \nu \, +\mu \, -\, 1) } E_{\alpha  \nu ,\alpha \,  \nu \, +\, \mu } (\pm   \lambda   t^{\alpha } ),} 
 \\  \phantom{\rule{1pt}{15pt}}  
 {{\ds \int _{0}^{\infty }}u^{\alpha } E_{\alpha ,\alpha \, +\, 1} (\pm   \lambda   u^{\alpha } )  \, f_{\nu ,\mu } (u)  \, du=
 \ t^{[(\alpha \, +\, 1)\,  \nu \, +\mu \, -\, 1] } E_{\alpha  \nu ,(\alpha \, +\, 1)\, \nu \, +\, \mu } (\pm   \lambda   t^{\alpha } ).} 
 \end{array} 
\end{equation} 
If $\beta = 1/2$, in the integrands of \eqref{GrindEQ__61_} the Mittag-Leffler functions are expressed then by the error functions [35]
\begin{equation} \label{GrindEQ__62_} 
\begin{array}{l} 
{E_{1/2} (\pm   z)=\e^{z^{2} } [1\pm \erf(z)]\quad ;\quad z=\lambda   u^{\alpha }, } 
\\  \phantom{\rule{1pt}{15pt}}  
{E_{1/2,1/2} (\pm   z)=\left\{\dfrac{1}{\sqrt{z} } \pm z \e^{z^{2} }
 [1\pm \erf(z)]\right\},}
  \\ \phantom{\rule{1pt}{15pt}}  
   {E_{1/2,3/2} (z)=\dfrac{ \e^{z} }{\sqrt{z} }\, \erf(z).}
    \end{array} 
\end{equation} 
If $\beta$  is positive integer, the Mittag-Leffler functions are expressed by elementary functions [36]. 
In particular cases, with 
$\nu=1/2,  \, \mu =0$. $\nu=\mu=1/2$ and $\nu=1/3,\, \mu=2/3$  
 the explicit form of the inverse transforms of exponential functions is known (see (15) and \eqref{GrindEQ__17_}).

 The Laplace transform of the error function is 
\begin{equation} \label{GrindEQ__63_} 
\begin{array}{l} {g(t)=
\erf \left(\dfrac{\lambda  }{ 2 t^{1/2} } \right)\quad ;\quad 0<\nu <1} 
\\  \phantom{\rule{1pt}{15pt}}  
{G(s)=\dfrac{1-\e^{-\, \lambda  s^{1/2} } }{s} \, \quad ;\quad G(s^{\nu } )
=
\dfrac{1-\e^{-\, \lambda  s^{\nu /2} } }{s^{\nu } } }
 \\ \phantom{\rule{1pt}{15pt}}  
  {\LL^{-1} \left\{\dfrac{1}{s^{\mu } } \bullet \dfrac{1-\e^{-\, \lambda  s^{\nu /2} } }{s^{\nu } } \right\}=  
  \LL^{-1} \, \left\{\dfrac{1}{s^{\mu } } -\dfrac{\lambda ^{2 \mu /\nu }
  \e^{-\, (\lambda ^{2/\nu }  s)^{\nu /2} } }{(\lambda ^{2/\nu } s)^{\mu } } \right\}=} 
  \\  \phantom{\rule{1pt}{15pt}}  
  {\dfrac{t^{\mu \, -\, 1} }{\Gamma (\mu )} -\lambda ^{2 (\mu \, -\, 1)/\nu } f_{\nu /2,\mu } \left(\dfrac{t^{\nu /2} }{\lambda ^{2/\nu } } \right),} 
  \end{array} 
\end{equation} 
which yields
\begin{equation} \label{GrindEQ__64_} 
{\ds \int _{0}^{\infty }}
\erf \left(\dfrac{\lambda  }{ 2 u^{1/2} } \right)  \, f_{\nu ,\mu } (u)  \, du
=
 \dfrac{t^{\mu \, -\, 1} }{\Gamma (\mu )} -\lambda ^{2 (\mu \, -\, 1)/\nu } f_{\nu /2,\mu } \left(\dfrac{t^{\nu /2} }{\lambda ^{2/\nu } }\right  ). 
\end{equation}


\noindent The Volterra functions are defined by the following integrals [36]
\begin{equation} \label{GrindEQ__65_} 
\begin{array}{l} 
{\nu (t)=
{\ds \int _{0}^{\infty }} \dfrac{t^{u} }{\Gamma (u+1)}  \, du,} 
\\ \phantom{\rule{1pt}{15pt}}  
 {\nu (t,\alpha )=
 {\ds \int _{0}^{\infty }}\dfrac{t^{u\, +\, \alpha } }{\Gamma (u+\alpha +1)}  \, du,} 
 \\ \phantom{\rule{1pt}{15pt}}  
 {\mu ((t,\beta ,\alpha )=
 {\ds \int _{0}^{\infty}}\dfrac{u^{\beta }  t^{u\, +\, \alpha } }{\Gamma (\beta +1) \Gamma (u+\alpha +1)}  \, du,} 
 \end{array} 
\end{equation} 
and their Laplace transforms are
\begin{equation} \label{GrindEQ__66_} 
\begin{array}{l} 
{L\left\{\nu (\lambda   t)\right\}=\dfrac{1}{s  \ln \left(\dfrac{s}{\lambda } \right)} \quad ;\quad \lambda >0} 
\\ \phantom{\rule{1pt}{15pt}}   
{\LL\left\{\nu (\lambda   t,\alpha )\right\}=\dfrac{\lambda ^{\alpha } }{s^{\alpha \, +\, 1}
   \ln \left(\dfrac{s}{\lambda } \right)} }
    \\ \phantom{\rule{1pt}{15pt}}  
     {\LL\left\{\mu ((\lambda   t,\beta ,\alpha )\right\}
     =\dfrac{\lambda ^{\alpha } }{s^{\alpha \, +\, 1}   \,
      \ln \left(\dfrac{s}{\lambda } \right)^{\beta \, +\, 1} }. } 
     \end{array} 
\end{equation} 
The logarithmic functions in the Laplace transforms permit to express the integrals of the Volterra functions with inverse Laplace transform of exponential function
in terms of convolution integrals.
 From
\begin{equation} \label{GrindEQ__67_} 
\begin{array}{l}
 {g(t)=\nu (\lambda   t)\quad ;\quad \lambda >0\quad ;\quad 0<\nu <1} 
 \\ \phantom{\rule{1pt}{15pt}}  
  {G(s)=\dfrac{1}{s  \ln \left(\dfrac{s}{\lambda } \right)} \quad ;\quad 
  G(s^{\nu } )=\dfrac{1}{s^{\nu }   \ln \left(\dfrac{s^{\nu } }{\lambda } \right)} }
   \\  \phantom{\rule{1pt}{15pt}}  
   { \LL^{-1} \left\{\dfrac{1}{s^{\mu } } \bullet \dfrac{1}{s^{\nu }   \ln [(s/\lambda ^{1/\nu } )^{\nu } ]} \right\}= 
   \LL^{-1} \left\{\dfrac{1}{s^{\nu \, +\, \mu \, -\, 1} } 
   \bullet \dfrac{1}{\nu  s  \ln (s/\lambda ^{1/\nu } )} \right\}=}
    \\ \phantom{\rule{1pt}{15pt}}  
     {\dfrac{t^{\nu \, +\, \mu \, -\, 2} }{\nu  \Gamma (\nu \, +\, \mu \, -1)} *\, \nu (\lambda ^{1/\nu } t)}
      \end{array} 
\end{equation} 
it follows that
\begin{equation} \label{GrindEQ__68_} 
{\ds \int _{0}^{\infty }}\nu (\lambda   u)  \, f_{\nu ,\mu } (u)  \, du
=
\dfrac{t^{\nu \, +\, \mu \, -\, 2} }{\nu  \Gamma (\nu \, +\, \mu \, -1)}
\, \star \, \nu (\lambda ^{1/\nu } t). 
\end{equation} 
Similarly from
\begin{equation} \label{GrindEQ__69_} 
\begin{array}{l}
 { g(t)=\nu (\lambda  t,\rho )\quad ;\quad \lambda ,\rho >0} 
 \\ \phantom{\rule{1pt}{15pt}}  
 {G(s)=\dfrac{\lambda ^{\rho } }{s ^{\rho \, +\, 1}  \ln \left(\dfrac{s}{\lambda } \right)} \quad ;\quad 
 G(s^{\nu } )=\dfrac{\lambda ^{\rho } }{s ^{(\rho \, +\, 1) \nu }  \ln \left(\dfrac{s^{\nu } }{\lambda } \right)} \quad ;\quad 0<\nu <1}
 \\ \phantom{\rule{1pt}{15pt}}  
 {\lambda ^{\rho } \LL^{-1} \left\{\dfrac{1}{s^{\mu } } \bullet \dfrac{1}{s ^{(\rho \, +\, 1) \nu }   \ln [(s/\lambda ^{1/\nu } )^{\nu } ]} \right\}
 = 
 \lambda ^{\rho } \LL^{-1} \left\{\dfrac{1}{s^{\mu \, -\, 1} } \bullet \dfrac{1}{\nu  s ^{(\rho \, +\, 1) \nu \, +\, 1}   \ln (s/\lambda ^{1/\nu } )} \right\}=} 
 \\ \phantom{\rule{1pt}{15pt}}  
  {\dfrac{t^{\mu \, -\, 2} }{\Gamma (\mu \, -\, 1)  \nu  \lambda } \, \star, \nu [\lambda ^{1/\nu } t,(\rho +1)  \nu ]}
   \end{array} 
\end{equation} 
we have
\begin{equation} \label{GrindEQ__70_} 
{\ds \int _{0}^{\infty }}\nu (\lambda  u,\rho )  \,f_{\nu, \mu } \, du
= \dfrac{t^{\mu \, -\, 2} }{\Gamma (\mu \, -\, 1)  \nu  \lambda } \, \star\, 
\nu [\lambda ^{1/\nu } t, (\rho +1)  \nu ]. 
\end{equation} 
Finally, the Laplace transform of the generalized Volterra function is
\begin{equation} \label{GrindEQ__71_} 
\begin{array}{l} 
{g(t)=\mu (\lambda  t,\xi ,\rho )\quad ;\quad \lambda ,\rho ,\xi >0\,;} 
\\
 {G(s)=\dfrac{\lambda ^{\rho } }{s ^{\rho \, +\, 1}  [\ln \left(\dfrac{s}{\lambda } \right)]^{\xi \, +\, 1} } \, ;\quad G(s^{\nu } )
 =\dfrac{\lambda ^{\rho } }{s ^{(\rho \, +\, 1) \nu }   [\ln \left(\dfrac{s^{\nu } }{\lambda } \right)]^{\xi \, +\, 1} } \, ;\quad 0<\nu <1, } 
 \\  \phantom{\rule{1pt}{15pt}}  
 {\lambda ^{\rho } \LL^{-1} \left\{\dfrac{1}{s^{\mu } } \bullet \dfrac{1}{s ^{(\rho \, +\, 1) \nu }   \{ \ln [(s/\lambda ^{1/\nu } )^{\nu } ]\} ^{\xi \, +\, 1} } \right\}= } 
 \\ \phantom{\rule{1pt}{15pt}}  
  {\dfrac{\lambda ^{\rho } }{\nu ^{\xi \, +\, 1} } \LL^{-1} \left\{\dfrac{1}{s^{\mu \, -\, 1} } \bullet \dfrac{1}{s ^{(\rho \, +\, 1) \nu \, +\, 1}   [\ln (s/\lambda ^{1/\nu } )]^{\xi \, +\, 1} } \right\}=}
   \\ \phantom{\rule{1pt}{15pt}}  
    {\dfrac{t^{\mu \, -\, 2} }{\Gamma (\mu \, -\, 1) \lambda  \nu ^{\xi \, +\, 1} } \, \star\, \mu [\lambda ^{1/\nu } t,\xi ,(\rho +1)  \nu ]} 
    \end{array} 
\end{equation} 
and \eqref{GrindEQ__65_} yields
\begin{equation} \label{GrindEQ__72_} 
{\ds \int _{0}^{\infty }}\mu (\lambda  u,\xi ,\rho )  f_{\nu, \mu }  \, du
=  \dfrac{t^{\mu \, -\, 2} }{\Gamma (\mu \, -\, 1) \lambda  \nu ^{\xi \, +\, 1} } \, \star\, \mu [\lambda ^{1/\nu } t,\xi ,(\rho +1)  \nu ] .
\end{equation} 

 The number of similar convolution integrals can be significantly enlarged if the Volterra functions are multiplied by $t^n$ with $n = 1, 2, 3, {\dots}$ , then their Laplace transforms should be differentiated $n$ times. 
 The result of such differentiations are linear combinations of these functions [33,36]. For example with $n = 1$ and $\lambda = 1$ we have
\begin{equation} \label{GrindEQ__73_} 
\begin{array}{l} 
{t  \nu (t,\rho )=
\LL^{-\, 1} \left\{-\dfrac{d}{d s} \left[\dfrac{1}{s^{\rho \, +\, 1}  \ln s} \right]\right\}
=\LL^{-\, 1} \left\{\dfrac{\rho +1}{s^{\rho \, +\, 2}  \ln s} +\dfrac{1}{s^{\rho \, +\, 2}  (\ln s)^{2} } \right\}=}
 \\ \phantom{\rule{1pt}{15pt}}  
  {(\rho +1)\, \nu [t,(\rho +1)  \nu ]+\, \mu [t,1,(\rho +1)  \nu ].}
   \end{array} 
\end{equation}

 \section*{Conclusions}

\noindent 

\noindent 

By applying the Efros theorem in the form established by Wlodarski it was possible to derive 
a number of  infinite integrals, finite integrals and integral identities with the function which represent the Laplace inverse transform of 
$s^{-\mu}\, \exp(-s^\nu)$ with $0<\nu<1$ and $0\le \mu < 1$. 
The extension to the cases $\mu \ge 1$ is dealt in Appendix

Derived by us integrals include in integrands elementary functions (power, exponential, logarithmic, trigonometric and hyperbolic functions) and the error functions, 
the Mittag-Leffler functions and the Volterra functions. 
Many results appear in form of the convolution integrals. 

Performing the inversion by the complex integration, it was possible to show that the inverse Laplace inverse transform, which means the original function,  can be also expressed in terms of the Wright functions and for particular values of parameters by the Mainardi functions. 

Using rules of operational calculus some properties of the inverse Laplace transform were derived.

\section*{Acknowledgments}
{The research work of  F.M.
has been carried out in the framework of the activities of the National Group of Mathematical Physics (GNFM, INdAM).
Both the authors would like to acknowledge
the unknown reviewers for their constrictive comments.}
\section*{Appendix: 
Additional properties of the inverse of the Laplace transform of
$s^{-\mu}\, \exp(-s^\nu)$
}

\noindent The restriction posed on the second parameter of the function 
$f_{\nu,\mu}(t)$, i.e. $\mu<1$ (in the inverse of the Laplace transform 
$s^{-\mu}\, \exp(-s^\nu)$)
 can be removed by using rules of operational calculus. 
 In the general case, for a non-negative  parameter 
 $\rho \ge 0$,  this inverse transform can always be expressed as the convolution integral which includes the corresponding power function:
$$
 \begin{array}{l}
  {f_{\nu ,\rho } (t) 
  =\LL^{-\, 1} \left\{\dfrac{\e^{-\,  s^{\nu } } }{s^{\rho } } \right\}
  =\LL^{-\, 1} \left\{\dfrac{1}{s^{\lambda } } \bullet \dfrac{\e^{- s^\nu } }{s^\mu } \right\}
  =\dfrac{t^{\lambda \, -\, 1} }{\Gamma (\lambda )} 
  \, \star \, f_{\nu ,\mu } (t)} \\ {\rho =\lambda +\mu } 
  \end{array}
  \eqno (A-1)
  $$
In the particular case $\mu  = 1$, it reduces to the simple integral,
 because $1/{s}$ expresses the integration operation:
$$
 f_{\nu ,1} (t)
 =\LL^{-\, 1} \left\{\dfrac{1}{s} \,\ e^{-\,  s^{\nu } } \right\}
 ={\ds \int _{0}^{t}}f_{\nu ,0} (u)\, du.
 \eqno(A-2)
$$
 The product of two identical or different Laplace inverse transforms can be derived in the form of convolution integrals from
$$
 \begin{array}{l} 
 {\LL^{-\, 1} \left\{\dfrac{\e^{-\,  s^{\nu } } }{s^{\mu } } \bullet \frac{\e^{-\,  s^{\nu } } }{s^{\mu } } \right\}
 =\LL^{-\, 1} \left\{\dfrac{\e^{-\,  2 s^{\nu } } }{s^{2 \mu } } \right\}
 =2^{2 \mu /\nu } \LL^{-\, 1} \left\{\dfrac{\e^{-\,  (2^{1/\nu }  s)^{\nu } } }{(2^{1/\nu } s)^{2 \mu } } \right\}} 
 \\  \phantom{\rule{1pt}{15pt}}
 {\LL^{-\, 1} \left\{\dfrac{\e^{-\,  s^{\nu } } }{s^{\mu } } \bullet \frac{\e^{-\,  s^{\nu } } }{s^{\rho } } \right\}
 =\LL^{-\, 1} \left\{\dfrac{\e^{-\,  2 s^{\nu } } }{s^{\mu \, +\, \rho } } \right\}=
 2^{(\mu \, \, +\, \rho )/\nu }
  \LL^{-\, 1} \left\{\dfrac{\e^{-\,  (2^{1/\nu }  s)^{\nu } } }{(2^{1/\nu } s)^{( \mu \, +\, \rho )} } \right\}}
   \end{array}
   \eqno (A-3)
$$
 which gives
$$
 \begin{array}{l} 
 {f_{\nu ,\mu } (t)\, \star \, _{\nu ,\mu } (t)
 ={\ds \int _{0}^{t}}f_{\nu ,\mu } (t-\xi )  f_{\nu ,\mu } (\xi ) d\xi
  =2^{(2 \mu \, -\, 1)/\nu } f_{\nu ,2\mu } \left(\frac{t}{2^{1/\nu } } \right)} 
 \\  \phantom{\rule{1pt}{15pt}}
  {f_{\nu ,\mu } (t)*f_{\nu ,\rho } (t)
  ={\ds \int _{0}^{t}}f_{\nu ,\mu } (t-\xi )  f_{\nu ,\rho } (\xi ) d\xi =2^{(\mu \, +\, \rho \, -\, 1)/\nu } f_{\nu ,\mu \, +\, \rho } \left(\frac{t}{2^{1/\nu } } \right)}
   \end{array}
   \eqno   (A-4)
$$
 These results can be generalized to \textit{n}-fold integrals when exists the factor 
 $1/s^{n}$, with ${n} = 1,2,3,{\dots}$.
  In the case of the product of identical Laplace inverse transforms it is reduced to the following convolution integral
$$
 \LL^{-\, 1} \left\{\frac{1}{s^{n} } F(s)\right\}=
 {\ds \int _{0}^{t}}\dfrac{(t-u)^{n\, -\, 1} }{(n-1)!}  \, f(u)\, du.
 \eqno (A-5)
 $$
 In similar way, if the inverse transform is in a more general form
$$
 \begin{array}{l} 
 {\LL\left\{f_{\nu ,\mu } (t);\alpha \right\}
 =
 \LL\left\{\dfrac{\e^{-\, \alpha  s^{\nu } } }{s^{\mu } } \right\}\quad ;\quad 
 \LL\left\{f_{\nu ,\rho } (t);\beta \right\}=
 \LL\left\{\dfrac{\e^{-\, \beta  s^{\nu } } }{s^{\mu } } \right\}\quad ;\quad \alpha ,\beta >0} 
 \\  \phantom{\rule{1pt}{15pt}}
  {\LL^{-\, 1} \left\{\dfrac{\e^{-\, \alpha  s^{\nu } } }{s^{\mu } } \bullet
   \dfrac{\e^{-\, \beta  s^{\nu } } }{s^{\rho } } \right\}
   =\LL^{-\, 1} \left\{\dfrac{\e^{-\, (\alpha \, +\, \beta ) s^{\nu } } }{s^{\mu \, +\, \rho } } \right\}=} 
  \\ \phantom{\rule{1pt}{15pt}}
   {(\alpha +\beta )^{(\mu \, +\rho )/\nu }\,
    \LL^{-\, 1} \left\{\dfrac{\e^{-[\, (\alpha \, +\, \beta )^{1/\nu }  s]^{\nu } } }
    {[(\alpha \, +\, \beta )^{1/\nu } s]^{\mu \, +\, \rho } } \right\}} 
    \end{array} 
    \eqno (A-6)
    $$
 we have
$$
 \begin{array}{l} 
 {f_{\nu ,\mu } (t;\alpha ) \,\star\, f_{\nu ,\rho } (t;\beta  )
 =
{\ds \int _{0}^{t}}f_{\nu ,\mu } (t-\xi ;\alpha )  f_{\nu ,\rho } (\xi ;\beta ) d\xi =}
  \\ \phantom{\rule{1pt}{15pt}}
   {(\alpha +\beta )^{(\mu \, +\, \rho \, -\, 1)/\nu } f_{\nu ,\mu \, +\, \rho } \left(\frac{t}{(\alpha +\beta )^{1/\nu } } ;(\alpha +\beta )\right)} 
   \end{array}
   \eqno (A-7)
$$
 The differentiation of Laplace inverse transforms with respect to the parameters 
$\nu$ and $\mu$
gives
$$
 \begin{array}{l} 
 {\LL\left\{\dfrac{\partial  f_{\nu ,\mu } (t)}{\partial  \nu } \right\}
 =
 -\dfrac{s^{\nu } \ e^{-\, s^{\nu } } \ln s}{s^{\mu } } =
 -\dfrac{ \e^{-\, s^{\nu } } }{s^{\mu \, -\, \nu \, -\, 1} } \bullet 
 \dfrac{\ln (s  C)}{s} +\frac{\ln C \e^{-\, s^{\nu } } }{s^{\mu \, -\, \nu \, } }\, ;
 \; C=\e^{\gamma } }
  \\   \phantom{\rule{1pt}{15pt}}
   {\LL\left\{\dfrac{\partial  f_{\nu ,\mu } (t)}{\partial  \mu } \right\}
   =-\dfrac{\mu   \e^{-\, s^{\nu } } }{s^{\mu \, +\, 1} } }
    \end{array}
    \eqno(A-8)
$$
 which yields the following inverses of (A-8)
$$
 \begin{array}{l}
  {\dfrac{\partial  f_{\nu ,\mu } (t)}{\partial  \nu } =
  \ln t \, \star \, f_{\nu ,\mu \, -\, \nu \, -\, 1} (t)+\gamma \, f_{\nu ,\mu \, -\, \nu } (t)} 
 \\  \phantom{\rule{1pt}{15pt}}
  {\LL\left\{\frac{\partial  f_{\nu ,\mu } (t)}{\partial  \mu } \right\}
  =-\mu \, f_{\nu ,\mu \, +\, 1} (t)
   =-\mu \, {\ds \int _{0}^{t}}f_{\nu ,\mu } (\xi )\, d\xi,  }
   \end{array}
   \eqno  (A-9)
$$
where $\gamma$  denotes the Euler constant.

\noindent 

\noindent

\section*{References}
\noindent \textbf{}

\noindent [1] Humbert, P. Nouvelle correspondances symboliques. 
Bull. Soc. Math., France,  69 (1945) 121--129.

\noindent [2] Pollard, H. 
The representation of $\exp(-x^\lambda)$ as a Laplace integral.
 Bull. Amer.  Math. Soc., 52 (1946) 908--910.

\noindent [3] Wlodarski, L. 
Sur une formule de Eftros. Studia Math., 13 (1952) 183--187.

\noindent [4] Mikusinski, J.
 Sur les fonctions exponentielles du calcul op\'{e}ratoire. 
 Studia Math.,  12 (1951) 208--224.

\noindent [5] Mikusinski, J.
 Sur la croissance de la function op\'{e}rationelle 
 $\exp(-s^\alpha \lambda)$. 
 Bull.  Acad. Polon. 4 (1953) 423--425.

\noindent [6] Mikusinski, J..
 Sur la function dont la transform\'{e}e de Laplace est 
 $\exp(-s^\alpha \lambda)$.
 Bull.  Acad. Polon. 6 (1958) 691--693..

\noindent [7] Mikusinski, J. 
On the function whose Laplace transform is $\exp(-s^\lambda)$. 
Studia Math.,  18 (1959) 195--198.

\noindent [8] Wintner, A.
 Cauchy's stable distributions and an ''explicit formula'' of Mellin.   
 Amer. J. Math., 78 (1956) 819--861.

\noindent [9] Ragab, F.M. 
The inverse Laplace transform of a exponential function. 
Comm.  Pure Appl. Math., 11 (1958) 115--127.

\noindent [10] Stankovi\v{c}, B. 
On the function of E.M. Wright. 
Publications de L'Institute  Math\'{e}matique. 
Nouvelle s\'{e}rie, 10 (1970) 113--124.  

\noindent [11] Mainardi, F.. 
Fractional Calculus and Waves in Linear Viscoelasticity, 
Imperial   College Press, London, 2010.

\noindent [12]  Hanyga, A. 
Multidimensional solutions of time-fractional diffusion-wave  equations. 
Proc. Royal Soc. A., London, 458 (2002) 933--957.

\noindent [13] Barkai, E. Fractional Fokker-Planck equation, solution, and application.
 Phys.  Rev. E. 63 (2001) 046118/1--17.  

\noindent [14] Penson, K.A., G\'{o}rska, K. 
Exact and explicit probability densities for one-sided  L\'{e}vy stable distribution. 
Phys. Rev. Letters, 105 (2010) 210604/1--9.

\noindent [15]  G\'{o}rska, K.,  Penson, K.A. 
L\'{e}vy stable distributions via associated integral  transform. 
J. Math. Phys., 53 (2012), 053302/1--10.   

\noindent [16] Efros, A.M. 
Some applications of operational calculus in analysis.
Mat. Sbornik,  42 (1935) 699--705.

\noindent [17]  Apelblat, A., Mainardi, F.
 Application of the Efros theorem to the Wright  functions of the second kind and other results. 
 Lecture Notes of TICMI 
(Tbilisi Int. Center of Mathematics and Informatics) , 
 21 (2020),  9--28.
Special issue edited by P.E. Ricci, dedicated to Prof. George Jaiani on the occasion of his 75th birthday, Tbilisi University Press, ISSN 1512-0511.

\noindent [18] Erd\'{e}lyi, A., Magnus, W., Oberhettinger, F., Tricomi, F.G. 
Tables of  Integral Transforms. McGraw-Hill, New York, 1954

\noindent [19] Roberts, G.E., Kaufman, H. 
Tables of Laplace Transforms. W.B. Saunders Co.,  Philadelphia, 1966.

\noindent [20] Hladik, J. La Transformation de Laplace a Plusieurs Variables.
 Masson et Cie  \'{E}diteurs. Paris, 1969.

\noindent [21] Oberhettinger, F., Badii, L. 
Tables of Laplace Transforms. Springer-Verlag,  Berlin, 1973.

 \noindent [22] Wright, E.M. On the coefficients of power series having exponential  singularities. J. London Math. Soc., 8 (1933) 71--79.

\noindent [23] Wright, E.M. 
The generalized Bessel function of order greater than one. 
Quart. J.  Math. Oxford, 11 (1940) 36--48.

\noindent [24] Montroll, E.W., Bendler, J.T. 
On L\'{e}vy (or stable) distributions and the Williams-Watts model of dielectric relaxation. 
J. Stat. Phys., 34 (1984) 129--192.

\noindent [25] Zolotarev, V.M. 
One-dimensional Stable Distributions. Translated from Russian  by H.H. McFaden. 
Amer. Math. Soc., Providence, Rhode Island, 1986.

\noindent [26] Mainardi, F.
On the initial value problem for the fractional diffusion-wave equation,
in S. Rionero and T. Ruggeri (Editors), 
 7th Conference on Waves and Stability in Continuous Media (WASCOM 1993),
World Scientific, Singapore, 1994, pp. 246--251. 


\noindent [27] Saichev, A.I., Zaslavsky, G.M. 
Fractional kinetic equations :solutions and  applications. 
Chaos, 7 (1997) 753--764.

\noindent [28] Mainardi, F., Tomirotti, M.. 
Seismitic  pulse propagation with constant Q and  stable probability distributions. 
Ann. Geofisica, 40 (1997) 1311--1328.

\noindent [29] Gorenflo, R., Luchko, Yu.,  Mainardi, F..
 Analytical properties and applications  of the Wright functions. 
 Fract. Calc. Appl. Anal., 2 (1999) 383--414.

\noindent [30] F. Mainardi, Yu. Luchko, G. Pagnini. 
The fundamental solution of the space-time  fractional diffusion equation. 
Fract. Calc. Appl. Anal., 4 (2001) 153--192. 

\noindent [31] Mainardi, F., Mura, A., Pagnini, G., Gorenflo, R. 
Fractional relaxation and time- fractional diffusion of distributed order. 
Proc. of the 2$^{nd}$ IFAC Workshop on  Fractional Differentiation and its Applications. 
Porto Portugal, 2006.

\noindent [33] Garg,, M., Rao, A.
 Fractional extensions of some boundary value problems in oil  strata. 
 Proc. India Acad. Sci. (Math. Sci.) 117 (2007) 267--281.

\noindent[33] Mainardi, F., Consiglio, A. 
The Wright functions of the second kind in Mathematical Physics.
Mathematics  (MDPI),  8 No 6 (2020),  884/1--26.
DOI: 10.3390/MATH8060884; E-print arXiv:2007.02098 [math.GM]

\noindent [34] Apelblat, A. 
Laplace Transforms and Their Applications. 
Nova Science  Publishers, Inc.,  New York, 2012.

\noindent [35] Apelblat, A., 
Differentiation of the Mittag-Leffler functions with respect to  parameters in the Laplace transform approach. 
Mathematics (MDPI), 8 (2020), 657/1--29.
DOI:10.3390/math8050657

\noindent [36] Gorenflo, R., Kilbas, A.A., Mainardi, F., Rogosin, S.V. Mittag-Leffler Functions,  Related Topics and Applications. $2^{nd}$ Edition, 
Springer,  Heidelberg, 2020.

\noindent [37] Apelblat, A. Volterra Functions. 
Nova Science Publishers, Inc., New York, 2008.

\noindent

\end{document}